\documentclass[11pt]{article}
\usepackage[margin=1.278in]{geometry}

\usepackage{amsmath,amsthm,amssymb,picins,graphicx}
\usepackage{hyperref}

\hypersetup{
pdfauthor = {Boris Bukh},
pdftitle = {Set families with a forbidden subposet},
pdfsubject = {Mathematics},
pdfkeywords = {Sperner's theorem, forbidden subposet, Boolean lattice, induced subposet, Hasse diagram, extremal set theory},
pdfstartview = {FitH},
pdfpagemode = {None}}
\newtheorem*{conjecture}{Conjecture}
\newtheorem{theorem}{Theorem}
\newtheorem{lemma}[theorem]{Lemma}
\newtheorem{proposition}[theorem]{Proposition}
\newcommand*{\F}{\mathcal{F}}
\newcommand*{\Z}{\mathbb{Z}}
\newcommand{\comment}[1]{}
\DeclareMathOperator{\pdist}{pdist}
\DeclareMathOperator{\Mon}{Mon}

\author{Boris Bukh
}
\title{Set families with a forbidden subposet}
\date{}

\newcommand*{\floor}[1]{\lfloor #1\rfloor}
\newcommand*{\abs}[1]{\lvert #1\rvert}
\DeclareMathOperator{\ex}{ex}

\begin{document}
\maketitle
\begin{abstract}
We asymptotically determine the size of the largest
family $\F$ of subsets of $\{1,\dotsc,n\}$
not containing a given poset $P$ if the Hasse
diagram of $P$ is a tree. This is a qualitative 
generalization of several known results including 
Sperner's theorem.
\end{abstract}
\section*{Introduction}
We say that a poset $P$ is a \emph{subposet} of a poset
$P'$ if there is an injective map $f\colon P\to P'$ such
that $a\leq_P b$ implies $f(a)\leq_{P'} f(b)$.
A poset $P$ is an \emph{induced subposet} of $P'$ if there is an injective
map $f\colon P\to P'$ for which $a\leq_P b$ if and only if
$f(a)\leq_{P'} f(b)$.  
For instance, \includegraphics[scale=0.32]{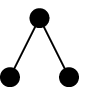} is a
subposet of \raisebox{-3pt}{\includegraphics[scale=0.32]{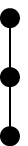}}, but not an induced
subposet. For a poset $P$ a \emph{Hasse diagram}, denoted by $H(P)$,
 is a graph whose vertices
are elements of $P$, and $xy$ is an edge if $x<y$ and for no other
element $z$ of $P$ we have $x<z<y$. 

Let $[n]=\{1,\dotsc,n\}$, and denote by $2^{[n]}$ the
collection of all subsets of $[n]$.  One
can think of a family $\F$ of subsets of $[n]$ as a poset
under inclusion. In this way $\F$ becomes an induced subposet of the Boolean
lattice.  In this paper we examine the size of the largest
family $\F\subset 2^{[n]}$ subject to the condition
that $\F$ does not contain a fixed finite subposet $P$. We do not
require $P$ to be an induced subposet.
A set family $\F$ not containing a subposet $P$ will be called a $P$-free family.
We denote by $\ex(P,n)$ the size of the largest
$P$-free family $\F\subset 2^{[n]}$.
For example, the classical
Sperner's theorem \cite{cite:sperner_orig} asserts that
$\ex(\includegraphics[scale=0.32]{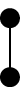}\,,n)=\binom{n}{\floor{n/2}}$.

Erd\H{o}s \cite{cite:erdos_littlewood_offord} extended Sperner's result, and
proved that if $C_l$ denotes the chain of length $l$, then $\ex(C_l,n)$ is
equal to the sum of $l-1$ largest binomial coefficients of order~$n$.
Katona and Tarj\'an\cite{cite:katona_tarjan} proved that
$\ex(\includegraphics[scale=0.32]{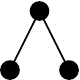}\,,n)=\binom{n}{\floor{n/2}}\bigl(1+O(1/n)\bigr)$. A common generalization of results of
Erd\H{o}s and Katona-Tarj\'an was established by
Thanh\cite{cite:thanh} who showed that if $P_{k,l}$ is any fixed poset with
vertex set $\{1,\dotsc,k\}\cup\{1',\dotsc,l'\}$ in
which the relations are $1<2<\dotsb<k$ and $1'<1$, $2'<1$, \dots
$l'<1$, then $\ex(P_{k,l})=
k\binom{n}{\floor{n/2}}\bigl(1+O(1/n)\bigr)$
(the error term was subsequently improved by
de Bonis and Katona\cite{cite:debonis_katona}). For example,
$\ex(\raisebox{-3pt}{\includegraphics[scale=0.32]{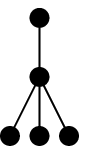}}\,,n)=2\binom{n}{\floor{n/2}}\bigl(1+O(1/n)\bigr)$. In \cite{cite:debonis_katona_swanepoel}
it is shown that $\ex(\includegraphics[scale=0.32]{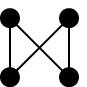}\,,n)=
\binom{n}{\floor{n/2}}+\binom{n}{\floor{n/2}+1}$.
 Griggs and Katona \cite{cite:griggs_katona} proved that $\ex(\includegraphics[scale=0.32]{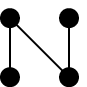}\,,n)=\binom{n}{\floor{n/2}}\bigl(1+O(1/n)\bigr)$. Recently, Griggs and Lu\cite{cite:griggs_lu} proved
that $\ex(L_{4k},n)=\binom{n}{\floor{n/2}}\bigl(1+O(1/n)\bigr)$, where $L_{4k}$ is the loop
of length $4k$ on two adjacent levels of the Boolean lattice.

Let $\Mon(\Z)$ be the set of all
functions $f\colon \Z\to\{0,1\}$ such that $f(n)=1$
and $f(-n)=0$ for all sufficiently large $n$. The
elements of $\Mon(\Z)$ will be called eventually monotone functions.
For $f,g\in\Mon(\Z)$ write $f\leq g$ if $f(n)\leq g(n)$ for all $n$.
Then $(\Mon(\Z),<)$ is a distributive lattice.
Note that
$\sum_n f(n)-g(n)$ is well-defined for every $f,g\in\Mon(\Z)$.
We define a \emph{level} of $\Mon(\Z)$ to be a maximal family 
$L\subset \Mon(\Z)$
satisfying $\sum_{n} f(n)-g(n)=0$ for all $f,g\in L$. Note that 
a level of $\Mon(\Z)$ is an antichain. The poset $\Mon(\Z)$
can be thought of as the induced subposet of $2^{\Z}$ spanned by the
set $\{X\subset \Z : \abs{X\triangle Y}<\infty\}$, 
for some fixed $Y\subset \Z$ which is neither finite nor cofinite.

The simplest explanation for all the results above is
the following conjecture.
\begin{conjecture}
For a finite poset $P$ let $l(P)$ be the maximum
number of levels in $\Mon(\Z)$ so that their
union does not contain $P$ as a subposet, then
\begin{equation*}
\ex(P,n)=l(P)\binom{n}{\floor{n/2}}\bigl(1+O(1/n)\bigr).
\end{equation*}
\end{conjecture}
Intuitively the conjecture asserts that the largest $P$-free
family is essentially the union of the maximum number of middle levels 
of the Boolean lattice $2^{[n]}$ not containing $P$. 
If true, the conjecture would be an analogue of Erd\H{o}s-Stone-Simonovits 
theorem from the extremal graph theory which asserts that
the largest graph not containing a given graph $G$ is
essentially the largest complete partite graph not containing $G$

In this paper, we establish the conjecture whenever $H(P)$ is a tree, 
generalizing several of the results mentioned above.
Unlike the papers above we are not concerned with establishing
the best possible bounds inside the $O(1/n)$ term, which allows us to give
a rather short proof.
The rest of the paper is occupied by
the proof of the following theorem.
\begin{theorem}\label{mainthm}
If $P$ is a finite poset and $H(P)$ is a tree, then
\begin{equation*}
\ex(P,n)=(h(P)-1)\binom{n}{\floor{n/2}}\bigl(1+O(1/n)\bigr)
\end{equation*}
where $h(P)$ is the height of $P$, i.e., the length of the longest
chain in $P$. Moreover, $l(P)=h(P)-1$ for such a~$P$.
\end{theorem}
For the case $h(P)=2$, the theorem was also independently proved by Griggs and Lu\cite{cite:griggs_lu} by a very
different argument.

\section*{Proof idea}
Before embarking on the proof of Theorem~\ref{mainthm}, we first non-rigorously sketch a simple proof for the special 
case $h(P)=2$. The proof unfortunately does not generalize to $h(P)\geq 3$, but it will motivate the otherwise
hard-to-follow technical details of the more general proof. We shall need a strengthening
of Sperner's lemma, to the effect that if $\abs{\F}>\binom{n}{\floor{n/2}}$, then not only
there are pairs of comparable sets, but a plentitude of such pairs.
The following statement is easy and can be proved similarly to Lemma~\ref{spernerlem} below
(see \cite{cite:kleitman} for a sharper result for small $\epsilon$).
\begin{lemma}\label{heighttwolem}
If $\abs{\F}\geq (1+\epsilon)\binom{n}{\floor{n/2}}$, then there are at least $(1/10)n \epsilon\abs{\F}$ pairs
of sets $F_1\subset F_2$ contained in $\F$. 
\end{lemma}
\begin{proof}[Proof of the case $h(P)=2$ of Theorem~\ref{mainthm} assuming Lemma~\ref{heighttwolem}]
Suppose $\F\subset 2^{[n]}$ is a family of $\abs{\F}\geq (1+20\abs{P}/n)\binom{n}{\floor{n/2}}$
sets, where $P$ is a poset. We will show that
$\F$ contains a copy of $P$. Let $G$ be the graph with vertex set 
$\F$ where a pair of distinct sets $F_1,F_2\in\F$ connected by an edge if $F_1\subset F_2$ or $F_2\subset F_1$. 
By Lemma~\ref{heighttwolem} the average degree of $G$ is at least $4\abs{P}$. 
Since every graph of
of average degree $d$ contains a non-empty subgraph of minimum degree at least $d/2$,
the graph $G$ contains a subgraph $G'$ of minimum degree at least $2\abs{P}$. 

Then we shall embed elements of $P$ one-by-one into $V(G')$, in such a way that at each step embedding is injective,
and preserves order. More precisely, we assume that $P$ is not a subposet of $V(G')$,
and use this assumption to construct a sequence of embeddings $\pi_i\colon P_i\to V(G')$, where
\renewcommand*{\theenumi}{(\alph{enumi})}
\renewcommand*{\labelenumi}{\alph{enumi})}
\begin{enumerate}
\item\label{partembedfirst} $P_i$ is an induced subposet of $P$, and $H(P_i)$ is a tree,
\item $P_i\setminus P_{i-1}$ consists of a single element,
\item\label{partembedlast} $v\leq_{P_i} u$ implies $\pi_i(v)\leq \pi_i(u)$ for all $u,v\in P_i$.
\end{enumerate} 

First, we embed some element of $P$ arbitrarily into $V(G')$, and let $P_1$ to consist
of that single element. Then, for 
each $i\geq 2$ let $v\in P\setminus P_{i-1}$ be a not-yet-embedded element of $P$,
which is comparable to some $u\in P_{i-1}$, and let $P_i=P_{i-1}\cup\{v\}$. Since $h(P)=2$, and $H(P)$ is a tree,
$u$ is the only element of $P_{i-1}$ comparable to $v$. 
We shall define embedding $\pi_i$ that agrees with $\pi_{i-1}$ on $P_{i-1}\setminus \{u\}$. 

Without loss of generality we can assume that $u<v$. Set $\tau_1=\pi_{i-1}$, and write 
$u_1=\tau_1(u)$. Since $u_1\in V(G')$, and degree of every vertex in $V(G')$ is at 
least $2\abs{P}\geq \abs{P}+1$ there is at least one 
neighbor $u_2$ of $u_1$ in $G'$ which is not in the image of $\tau_1$. If $u_1<u_2$, then 
we let $\pi_i$ to be an extension of $\tau_1$ sending $v$ to $u_2$.  
Suppose $u_1>u_2$. Let $\tau_2$ be the embedding obtained
from $\tau_1$ by mapping $u$ to $u_2$ instead of $u_1$. The embedding $\tau_2$ satisfies the same properties
\ref{partembedfirst}--\ref{partembedlast} that $\pi_{i-1}$ does. As $2\abs{P}\geq \abs{P}+2$ we 
can again find a neighbor $u_3$ of $u_2$ which is neither $u_1$ nor in the image of $\pi_2$,
and we can assume that $u_2>u_3$. Repeating this process $\abs{P}$ times yields a chain $u_1>u_2>\dotsc>u_{\abs{P}}$
of elements of $V(G')$.  However $P$ is a subposet (but not necessarily an induced subposet) 
of any linear extension of itself, and thus embeds into
$\{u_1,\dotsc,u_{\abs{P}}\}\subset \F$.
\end{proof}
The proof above is clearly similar to the proof that 
every tree on $d$ vertices embeds into a graph of average 
degree $2d$. To extend the
proof above to the case $h(P)=3$, for example, it is tempting to 
replace the graph $G$ by a $3$-uniform hypergraph of triples of 
sets in a chain. The problem with this approach is lack of any good replacement for 
the concept of minimum degree. The solution is therefore to
eliminate minimum degree from the proof 
entirely. To see how it is done, we present an alternative way of 
embedding trees in graphs of large
average degree. It is far more wasteful, but avoids minimum degrees.
\begin{proposition}\label{treeembedprop}
For every tree $T$ there is a $d_0(T)$ such that $T$ embeds into every finite graph of average degree at least~$d_0(T)$. 
\end{proposition}
\begin{proof}
The proof is by induction on $T$, with $\abs{T}=1$ being the trivial base case. Assume $\abs{T}\geq 2$. Let $v$ be any leaf of $T$,
and $u$ be its unique neighbor. Set $T'=T\setminus v$. We will show that we can take $d_0(T)=d_0(T')+2\abs{T}$. Suppose
$G$ is of average degree at least $d_0(T')+2\abs{T}$. Let $B=\{v \in G : \deg_G(v)\leq \abs{T} \}$. Then the graph $G'=G\setminus B$
is only $\abs{B}\abs{T}$ edges smaller than $G$, and thus has has average degree at least 
$(d_0(T')+2\abs{T})-2\abs{B}\abs{T}/\abs{V}\geq d_0(T')$. By the induction hypothesis there is an embedding $\pi\colon T'\to G'$.
As $\pi(u)\in G'$, we infer $\deg_G(\pi(u))>\abs{T}$, implying that there is at least one neighbor of $\pi(u)$ that is not in $\pi(T')$. Use this neighbor to extend the embedding $\pi$ of $T'$ to an embedding of $T'\cup\{v\}=T$.
\end{proof}
For technical reasons, it turns out to be easier to work
with marked chains rather than hypergraphs; intuitively this change corresponds to
hypergraphs with weighted edges.  
In the other respects, the proof of Theorem~\ref{mainthm} is a straightforward generalization of the 
two arguments above. 

\section*{Proof of Theorem~\protect\ref{mainthm}}
By an \emph{interval} in a poset $P$ we mean a
set of the form $[x,y]=\{z\in P : x\leq z\leq y\}$. 
A \emph{maximal chain} in a poset $P$ is a chain, which
is not contained in any other chain. In particular, 
a maximal chain in $2^{[n]}$ is a chain of sets
$\emptyset=S_0\subset S_1\subset \dotsb \subset S_n=[n]$
with $\abs{S_i}=i$. A \emph{$k$-marked chain} with
markers $F_1,\dotsc,F_k$ is a $k+1$-tuple
$(M,F_1,\dotsc,F_k)$ where $M$ is a maximal chain in $2^{[n]}$,
$F_1\supset\dotsb\supset F_k$ and $F_1,\dotsc,F_k$ belong to~$M$.
\begin{lemma}\label{spernerlem}
If $\F\subset 2^{[n]}$ is of size
\begin{equation*}
\abs{\F}\geq (k-1+\epsilon)\binom{n}{\floor{n/2}},
\end{equation*}
then there are at least $(\epsilon/k)n!$
 $k$-marked chains whose markers belong to $\F$.
\end{lemma}
\begin{proof}
Let $C_i$ be the number of maximal chains
that contain exactly $i$ sets from $\F$.
Counting the number of pairs $(F,M)$ where $F\in \F$ is an element of a maximal chain $M$ in
two different ways, we obtain
\begin{equation*}
\sum_i i C_i=\sum_{F\in\F}\frac{n!}{\binom{n}{\abs{F}}}\geq \abs{\F}\frac{n!}{\binom{n}{\floor{n/2}}} \geq (k-1+\epsilon)n!.
\end{equation*}
From this and $\sum C_i=n!$, we infer that
\begin{equation*}
\sum_{i\geq k} i C_i\geq \sum_i i C_i -(k-1)\sum_{i\leq k-1} C_i 
\geq \epsilon n!.
\end{equation*}
The number of $k$-marked chains with markers in $\F$ is
\begin{align*}
\sum_{i\geq k}\binom{i}{k}C_i&=\sum_{i\geq k}\binom{i-1}{k-1}\frac{i}{k}C_i
\geq \frac{1}{k}\sum_{i\geq k}i C_i\geq \frac{\epsilon}{k}n!.\qedhere
\end{align*}
\end{proof}
Call a poset $P$ of height $k$ \emph{saturated} if every maximal chain
is of length $k$. For us the maximal chains in $P$ play the role
analogous to the edges of a tree $T$ in Proposition~\ref{treeembedprop}.
However, in general the edges might have different sizes, which is
analogous to dealing with non-uniform hypergraphs. The saturated posets are 
the analogues of uniform hypergraphs.

The next two lemmas establish a couple of intuitively obvious, but annoyingly hard to 
rigorously prove facts
about saturated posets whose Hasse diagram is a tree. The first lemma will be used
to reduce the problem of embedding an arbitrary $P$ to the problem of embedding
saturated~$P$. The second lemma will allow us to do induction on $\abs{P}$.
\begin{lemma}\label{saturcont}
If $P$ is a finite poset with $H(P)$ being a tree, then $P$ is an induced subposet of some saturated
finite poset $\tilde{P}$ with $H(\tilde{P})$ being a tree, and $h(P)=h(\tilde{P})$.
\end{lemma}
\begin{proof}
For the purpose of this proof let $s(P)$ be the number of maximal 
chains in $P$. Since every element is contained in some maximal chain,
$\abs{P}\leq s(P)h(P)$, implying that for fixed $s$ and $k$
there only finitely many posets $P$ with $s(P)=s$ and $h(P)=k$.
Assume there is a counterexample to the lemma. Let
$P$ be a counterexample with the largest number of elements
for given $s(P)$ and $h(P)$.

Since $P$ is not saturated there is a maximal 
chain $v_1<\dotsb<v_t$ in $P$ of length $t<k$.
For each $i=0,\dotsc,t-1$ we can define a new poset $P_i$
which is obtained from $P$ by adding a new element
$v$ and two new relations $v_i<v$ and $v<v_{i+1}$
(in the case $i=0$ we add only one new relation).
Clearly, $s(P_i)=s(P)$ and $P$ is an induced subposet of $P_i$.
We will show by induction on $i$ that for each
$i=0,\dotsc,t-1$ either $h(P_j)=h(P)$ for some $j\leq i$,
or there is a chain in $P$ of length $k-i$ whose smallest element
is $v_{i+1}$.

If $h(P_0)>k$, then there is a chain $C$ in $P_0$ of length $k+1$.
Since $C$ is not a chain of $P$, it contains $v$. Thus $v_1\in C$, and 
$C\setminus \{v\}$ contains $v_1$ and has length $k$.
Now suppose $i\geq 1$ and we have established the inductive claim for
all smaller values of $i$. If $h(P_j)>h(P)$ for all $j=0,\dotsc,i-1$, then
by the induction hypothesis there is a chain
$C\subset P$ of length $k-i+1$ whose smallest element is $v_i$. Therefore,
all chains in $P$ whose largest element is $v_i$ have lengths
not exceeding $i$. Since the length of the chain $v_1<\dotsb<v_i$ is $i$, 
the assumption $h(P_i)>k$ implies that there is a chain in $P$ of length at
least $k-i$ whose smallest element is $v_{i+1}$. This establishes the
inductive claim.

Hence, if $h(P_i)>h(P)$ for all $i$, there is a chain of length $k-t+1\geq 2$
whose smallest element is $v_t$. This contradicts the maximality
of $v_1<\dotsb<v_t$, implying that $h(P_i)=h(P)$ for some $i$.
However, $P$ was assumed to be the largest counterexample with given $s(P)$ and $h(P)$.
Therefore, there are no counterexamples to the lemma.
\end{proof}
\begin{lemma}\label{leaflemma}
Suppose $P$ is a saturated finite poset of height $h(P)=k\geq 2$, and $H(P)$ is a tree.
Furthermore, assume that $P$ is not a chain.
Then there is a $v\in P$, which is a leaf in $H(P)$,
and an interval $I$ of length $\abs{I}\leq k-1$ 
containing $v$ such that $H(P\setminus I)$
is a tree, and the poset
$P\setminus I$ is a saturated poset of height $h(P)$.
\end{lemma}
\begin{proof}
A sequence $v_1,\dotsc,v_l\in P$ is said to be a \emph{poset path}
of length $l-1$ if $v_i$ and $v_{i+1}$ are comparable for all $i$. 
A \emph{poset distance} between $v$ and $u$, denoted $\pdist(v,u)$,
is the shortest length of a poset path connecting them.

\parpic[r]{\includegraphics{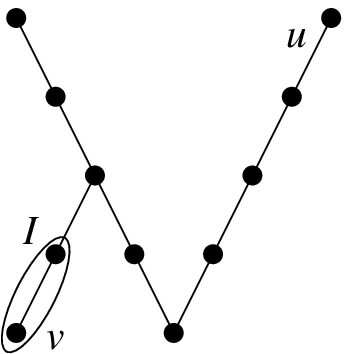}}
Let $v$ and $u$ be a pair of leaves maximizing $\pdist(v,u)$,
and let $v=v_1,v_2,v_3,\dotsc,v_l=u$ be the shortest
poset path between them. Observe that $\pdist(v,u)\geq 2$, 
for $P$ is not a chain. Without loss of generality
we can assume that $v<v_2$, for we can consider the opposite
poset otherwise. Let $I_0$ be the 
longest interval containing $v$, all of whose
elements have degree at most two in $H(P)$. 
If $\abs{I_0}< h(P)$, set $I=I_0$. If $\abs{I_0}=h(P)$,
set $I=I_0\setminus\{\max I_0\}$, where $\max I_0$ is the
largest element of $I_0$. 

Let $C$ be an arbitrary maximal chain in $P\setminus I$. We will show
that $C$ has length $h(P)$. 
Suppose $C$ contains an element $w$ which is incomparable
with $I$. Let then $M$ be a maximal chain in $P$ containing
$C$. Since $w$ is comparable to every element of $M$, the chain
$M$ is disjoint from $I$, implying $C=M$ and $h(C)=h(M)=h(P)$.
So we can suppose that all elements of $C$ are comparable to $I$.

If $\abs{I_0}=h(P)$, then the only element in $P\setminus I$
 which is comparable to $I$ is $\max I_0$. Since $P$ is not 
a chain, there is a chain of length two containing $\max I_0$,
which is disjoint from $I$, contradicting the maximality of $C$. 
Hence, we can assume that $\abs{I_0}<h(P)$. 
Let $w$ be any element of $P\setminus I$ comparable to an element $z$
of $I$. If $w<z$, then $\min([w,z]\cap I)$ has degree at least $3$ in $H(P)$.
If $z<w$ and $\max I\not<w$, then $\max([z,w]\cap I)$ has degree at least $3$ in $H(P)$.
In either case, we reach a contradiction with the choice of~$I$.
Hence $w$ exceeds all elements of $I$. Taking $w=\min C$, it follows
that $\max I<\min C$. By the maximality of $C$, there is no $z\in P\setminus I$ satisfying
 $\max I<z<\min C$. Upon taking $w=v_2$, it also follows that $\max I<v_2$, implying
$\min C\leq v_2$. If $v_2=\min C$, then $C\cup\{v_3\}$ is also a chain, 
contradicting the maximality of $C$. Let $y=\max C$. The only path from $u$ to $y$
in $H(P)$ goes through $v_2$ and $\min C$. Therefore every poset path
from $u$ to $y$ has to go through $v_2$ and $\min C$, implying $\pdist(y,u)>\pdist(v,u)$.
Though the element $y$ needs not to be a leaf itself, if $z$ is
any leaf of $H(P)$ such that the path from $u$ to $z$ goes through $y$,
then $\pdist(z,u)\geq \pdist(y,u)>\pdist(v,u)$, contradicting the choice of $v$ and $u$.
\end{proof}
The core of the proof of Theorem~\ref{mainthm} is contained in the
following lemma. In the lemma the family of $k$-marked chains $\mathcal{L}$
plays analogous role to the graph $G$ in the Proposition~\ref{treeembedprop},
with $\F$ being analogous to the vertex set of~$G$. The condition that 
$\F$ does not contain a chain of length $K$ comes from the fact that every 
poset embeds into every sufficiently long chain.
\begin{lemma}\label{mainlemma}
Let $P$ be a saturated finite poset of height $h(P)=k\geq 2$,  
whose Hasse diagram is a tree. Suppose $\F\subset 2^{[n]}$ is a set family, 
such that no chain
contains more than $K$ sets from $\F$, and all sets in $\F$
are of size between $n/4$ and $3n/4$. Moreover, suppose
$\mathcal{L}$ is a family of $k$-marked chains with markers in $\F$
of size
\begin{equation*}
\abs{\mathcal{L}}\geq \frac{\binom{\abs{P}+1}{2}4^{K+1}}{n} n!.
\end{equation*}
Then there is an embedding of $P$ into $\F$ in which every maximal chain
of $P$ is mapped to the set of markers of some $k$-marked chain in $\mathcal{L}$.
\end{lemma}
\begin{proof}
The proof is by induction on $\abs{P}$. If $P$ is the chain of length $k$,
then finding the required embedding is easy: marked elements on any $L\in\mathcal{L}$
form the desired chain. Now suppose we want to embed $P$, and have already
established the lemma for all smaller saturated posets. Use the preceding lemma to obtain
a leaf $v$ and an interval $I\ni v$ such that $P\setminus I$ 
is a still a saturated poset of
height $k$. By passing to the opposite poset to $P$
and replacing $\F$ by $\bar{\F}=\{[n]\setminus F: F\in\F\}$
if necessary, we can assume that $v$ is smaller than any 
element that is comparable with $v$.
Let $C$ be a maximal chain containing $I$. Let 
$s=\abs{C\setminus I}=k-\abs{I}$. Note that $s\geq 1$.

Call a chain $F_1\supset \dotsb \supset F_s$ of length $s$
a \emph{bottleneck} if there is a set $\mathcal{S}\subset \F$
with than $\abs{P}$ elements such that for every 
$k$-marked chain in $\mathcal{L}$ of the form 
$(M,F_1,\dotsc,F_s,F_{s+1},\dotsc,F_k)$ we have
$\mathcal{S}\cap\{F_{s+1},\dotsc,F_k\}\neq\emptyset$. Such an $\mathcal{S}$
is said to be a \emph{witness} to the fact that $F_1\supset\dotsb\supset F_s$
is a bottleneck. Note that without loss of generality, a witness contains 
only proper subsets of $F_s$. For each bottleneck $F_1\supset\dotsb\supset F_s$, let $\mathcal{S}(F_1,\dotsc,F_s)$
be a fixed witness containing only proper subsets of $F_s$. 
Call a $k$-marked chain $(M,F_1,\dotsc,F_k)\in\mathcal{L}$ \emph{bad} 
if for some $s$ the chain $F_1\supset\dotsb\supset F_s$ is a bottleneck.

Consider any $s$-element set $R=\{r_1,\dotsc,r_s\}$
of integers with $1\leq r_1<\dotsb<r_s\leq K$. If $M$ is a maximal chain in $2^{[n]}$
containing at least $r_s$ elements from $\F$, let 
$F_1\supset F_2\supset \dotsb\supset F_{r_s}$ be the $r_s$ largest of these elements.
The subchain of $F_1\supset F_2\supset \dotsb\supset F_{r_s}$ indexed by $R$ is
$F_{r_1}\supset F_{r_2}\supset\dotsb\supset F_{r_s}$, and we denote it by $C_{R}(M)$.
If $C_{R}(M)$ is a bottleneck, and $\mathcal{L}$ contains a
$k$-marked chain of the form $(M,\dotsc)$, whose $s$ largest markers
are $C_R(M)$, then we say that $M$ is \emph{$R$-bad}. Intuitively, the $R$-bad
chains correspond to the edges adjacent to the vertices of low degree 
in Proposition~\ref{treeembedprop}.  

Pick a maximal chain $M$ in $2^{[n]}$ uniformly at random.
Let $B_{R}$ be the event that $M$ is $R$-bad. We will estimate 
$\Pr[B_{R}]$ for each fixed $R$ individually.

One way to pick a random maximal chain $M$ of $2^{[n]}$ is to start with $[n]$
and remove elements one by one, each step choosing an element uniformly 
at random among the remaining elements. Thus one can generate chain $M$ in two stages.
In the first stage, we remove elements from $[n]$ at random until either we encounter
$r_s$ sets from $\F$, or until we run out of elements to remove. Denote by $T$ the chain obtained at
the end of the first stage (it is not a maximal chain, unless we ran out of elements). 
In the second stage, we resume removing elements at random from $\min T$, until
no elements are left. If $T$ is not maximal, then $C_{R}(M)$ is independent of what happens
in the second stage, and $C_{R}(T)$ is defined in the obvious way.

If $T$ is a maximal chain, or $C_R(T)$ is not a bottleneck, then $B_R$ does not hold.
Otherwise, let  $\mathcal{S}=\mathcal{S}(C_R(T))$ be the witness that $C_R(T)$ is a bottleneck. Recall that 
$\mathcal{S}\subset \F$,
$\abs{\mathcal{S}}<\abs{P}$ and $\mathcal{S}$ meets every $k$-chain in $\mathcal{L}$ whose top $s$
markers are $C_R(T)$. Let 
\begin{equation*}
\mathcal{T}_R=\{\text{chain }T_0\text{ in }2^{[n]} : \abs{T_0\cap \F}=r_s, C_R(T_0)\text{ is a bottleneck}\}.
\end{equation*}
The probability that $M$ meets $S$ is thus
\begin{align*}
\Pr[M\cap \mathcal{S}\neq\emptyset\, |\,  C_{R}(T)\text{ is a bottleneck} ]&\leq 
\max_{T_0\in\mathcal{T}_R}\Pr[M\cap \mathcal{S}\neq\emptyset | T=T_0]\\
&\leq \max_{T_0\in\mathcal{T}_R}\abs{\mathcal{S}(C_R(T_0))}\max_{\substack{F\in\F\\F<\min T_0}}\Pr[F\in M | T=T_0]\\
&\leq \abs{P}\max_{T_0\in\mathcal{T}_R}\max_{\substack{F\in\F\\F<\min T_0}}\Pr[F\in M | T=T_0]\\
&\leq \abs{P}\max_{T\in\mathcal{T}_R}\max_{\substack{F\in\F\\F<\min T_0}} \frac{1}{\abs{F}+1}\\
&\leq \abs{P}\max_{F\in\mathcal{F}} \frac{1}{\abs{F}+1},\\
&\leq 4\abs{P}/n
\end{align*}
where the fourth inequality follows because at the step before 
obtaining $F$, we have $\abs{F}+1$ choices as to which element 
to remove, with at most one choice yielding $F$.
If $M\cap \mathcal{S}=\emptyset$, then there is no $k$-marked chain of the form  $(M,\dotsc)$ in $\mathcal{L}$,
and $B_{R}$ does not hold. Therefore
\begin{align*}
\Pr[B_{R}]&=\Pr[C_{R}(T)\text{ is a bottleneck}]\Pr[B_R|C_{R}(T)\text{ is a bottleneck}]\\
&\leq \Pr[C_{R}(T)\text{ is a bottleneck}]\Pr[M\cap \mathcal{S}\neq\emptyset\, |\, C_{R}(T)\text{ is a bottleneck}]\leq 4\abs{P}/n.
\end{align*}
Since $R$ is a subset of $[K]$, the number of pairs $(M,R)$ where
$M$ is an $R$-bad maximal chain is at most $(4\abs{P}\binom{K}{s}/n)n!$.
Since no chain contains more than $K$ elements of $\F$, 
every bad $k$-marked chain gives rise to one
such pair $(M,R)$. Since $R\subset[K]$,
every pair $(M,R)$ arises in at most $\binom{K}{s}$ ways, implying that
there are no more than $(4\abs{P}\binom{K}{s}^2/n)n!
\leq (\abs{P}4^{K+1}/n)n!$
bad $k$-marked chains in $\mathcal{L}$.

Let $\mathcal{L}'$ be the set of all good $k$-marked chains in $\mathcal{L}$.
There are 
\begin{equation*}
\abs{\mathcal{L}'}\geq \abs{\mathcal{L}}-\frac{\abs{P}4^{K+1}}{n}n!\geq
\frac{4^{K+1}\bigl[\binom{\abs{P}+1}{2}-\abs{P}\bigr]}{n}n!
=\frac{\binom{\abs{P}}{2}4^{K+1}}{n} n!
\end{equation*}
of them.
By the induction hypothesis there is an embedding 
$\pi\colon P\setminus I\to \F$. 
Recall that $C$ was a maximal chain containing $I$, and
look at $C\setminus I$. Since $P\setminus I$ is saturated, $C\setminus I$
is contained in some chain $C'$ of length $k$ in $P\setminus I$.
Therefore $\pi(C')$ is contained in some $L\in\mathcal{L}'$.
Since all $k$-marked chains in $\mathcal{L}'$ are good, $\pi(C\setminus I)$
is not a bottleneck. In particular, since $\abs{P\setminus I}<\abs{P}$, we infer that
$\pi(P\setminus I)$ is not a witness
that $\pi(C\setminus I)$ is a bottleneck. 
Thus there is a $k$-marked chain 
$\tilde{L}\in\mathcal{L}$ containing $\pi(C\setminus I)$ as markers, but not containing 
any other element of $\pi(P\setminus I)$ as a marker. Therefore, we
can map $I$ to the bottom $k-s$ markers of $\tilde{L}$, completing
the desired embedding.
\end{proof}
With most of the work already done, we are ready to prove our main result.
\begin{proof}[Proof of Theorem~\ref{mainthm}]
If $h(P)=1$, then $P$ is a single-element poset, and the theorem is trivially
true. So, assume that $h(P)\geq 2$. 
Consider the case when $P$ is a saturated poset, and suppose
\begin{equation*}
\abs{\F}\geq (h(P)-1)\binom{n}{\lfloor n/2\rfloor}\Bigl(1+\frac{h(P)\abs{P}^2 4^{\abs{P}+2}}{n}\Bigr)
\end{equation*}
and $n$ is sufficiently large. We will show that $\F$ contains $P$. The number of
sets $F\in 2^{[n]}$ with fewer than $n/4$ or more than $3n/4$ elements is
equal $2^n$ times the probability that for a randomly chosen $F\in 2^{[n]}$
we have $\abs{\abs{F}-n/2}>n/4$. Thus by Chernoff's
inequality the number of such sets $F\in 2^{[n]}$ 
is at most $2^n\cdot 2\exp\bigl(-2(n/4)^2/n\bigr)=o\bigl(\binom{n}{\lfloor n/2\rfloor}/n\bigr)$. Let
$\F'=\{F\in F: \abs{F-n/2}\leq n/4\}$. As our $n$ is sufficiently large,
\begin{equation*}
\abs{\F'}\geq (h(P)-1)\binom{n}{\lfloor n/2\rfloor}\Bigl(1+\frac{h(P)\abs{P}^2 4^{\abs{P}+1}}{n}\Bigr).
\end{equation*}
Therefore, from Lemma~\ref{spernerlem} and 
Lemma~\ref{mainlemma} applied with $k=h(P)$ and $K=\abs{P}$ it follows 
that  either there is an embedding of $P$ into $\F'$ or $\F'$ contains a chain
$C$ of length $\abs{P}$. In the latter case, we can find an embedding $\pi$
of $P$ into $\F'$ anyway by simply letting $\pi\colon P\to C$ be
any linear extension of $P$.

If $H(P)$ is a tree and $P$ is not a saturated poset, then by Lemma~\ref{saturcont} it is contained 
in some saturated poset $P'$ of height $h(P')=h(P)$, such that $H(P')$ is a tree. 
Therefore $\ex(P,n)=(h(P)-1)\binom{n}{\lfloor n/2\rfloor}(1+O(1/n))$
for every poset $P$, for which $H(P)$ is a tree.

It remains to prove that $l(P)=h(P)-1$. The inequality $l(P)\geq h(P)-1$
is clear, as a union of $h(P)-1$ levels does not contain a chain of length
exceeding $h(P)-1$, and hence does not contain $P$. Let $L_1,\dotsc,L_h$
be $h$ distinct levels of $\Mon(\Z)$, and let $L$ be their union.
Suppose furthermore that the levels $L_i$ are so ordered that 
for any functions $f_i\in L_i$, the inequality $\sum_n f_i(n)-f_j(n)>0$ holds
whenever $i>j$ (by the definition of a level, if the inequality holds 
between a pair functions in levels $L_i$ and $L_j$, then it holds for
all pairs).
 
To complete the proof, we need to exhibit
an embedding of $P$ into $L$. By Lemma~\ref{saturcont} it suffices to treat the case
when $P$ is saturated. We will prove the existence of the embedding by induction
on $\abs{P}$. If $P$ is a chain of length $h$, then the embedding is obvious.
Suppose $P$ is not a chain, we can find embedding for all smaller saturated $P$
of height $h$. By Lemma~\ref{leaflemma} there is a leaf $v$ and an interval
$I$ of the form $I=[v,u)$ such that $P\setminus I$ is a saturated poset of height $h$.
By induction $P\setminus I$ is embeddable into $L$. 
Fix any such embedding. Since $\pi(u)$ is contained in a chain of length
$k$ in $L$, and $P$ is saturated, it follows that $\pi(u)\in L_{\abs{I}+1}$.
Let $n_0$ be a large enough that $(\pi(w))(n)=1$ for all $w\in P\setminus I$ and 
$n\geq n_0$. Complete the embedding by mapping 
the interval $I$ to the interval of functions $f_1,\dotsc,f_{\abs{I}}\in\Mon(\Z)$
defined by 
\begin{equation*}
f_i(n)=\begin{cases}0,&\text{if }n_0\leq n\leq n_0+i-1,\\
(\pi(w))(n),&\text{otherwise}.\end{cases}\qedhere
\end{equation*}
\end{proof}

\section*{Concluding remarks}
Though it would be interesting to determine exactly or 
find very good asymptotic estimates for $\ex(P,n)$ in general, 
a first step is to find the leading term in the asymptotic.
In this paper we found the leading term of $\ex(P,n)$ whenever
$H(P)$ is a tree. For some posets $P$ whose Hasse diagram is not a tree,
one can find a poset $P'$ that contains $P$ and
 whose Hasse diagram is a tree with $l(P)=l(P')$,
to obtain that $\ex(P,n)\sim \ex(P',n)$. For example, 
$\ex(\includegraphics[scale=0.32]{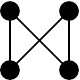}\,,n)\sim 
\ex(\raisebox{-3pt}{\includegraphics[scale=0.32]{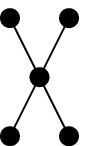}}\,,n)$,
and similarly for other complete two-level posets, thus recovering
the results from \cite[Section~5]{cite:debonis_katona}. 
The simplest two posets that are not subposets of trees with the same 
value of $l(P)$  are 
\raisebox{-3pt}{\includegraphics[scale=0.32]{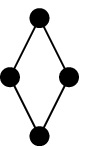}} and
\includegraphics[scale=0.32]{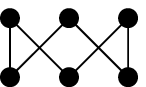}, and the asymptotics of the function $\ex$
for these posets is not known.

It is conceivable that the conjecture in this paper is even true
if its premise that $\F$ does not contain a subposet $P$ is replaced by
the weaker premise that $\F$ does not contain $P$ as an induced subposet.

\textbf{Acknowledgement. } I thank M\'at\'e Matolcsi for reading a preliminary 
version of this paper, and two referees for useful suggestions.

\bibliographystyle{alpha}
\bibliography{treeposetturan}

\end{document}